\documentclass[final,12pt]{article}
\usepackage{amsfonts,amssymb,amsthm,amsmath,mathtools}
\usepackage{pslatex,a4wide}
\usepackage{booktabs,float,graphicx}
\usepackage{tikz}
\usepackage[noadjust]{cite}
\usepackage{hyperref}

\newtheorem{theorem}{Theorem}[section]
\newtheorem{lemma}[theorem]{Lemma}
\newtheorem{corollary}[theorem]{Corollary}

\theoremstyle{definition}
\newtheorem{conjecture}[theorem]{Conjecture}

\theoremstyle{remark}

\numberwithin{equation}{section}

\newcommand{\diag}{\operatorname{diag}}
\newcommand{\bd}{\mathbf{d}}
\newcommand{\T}{\mathsf{T}}
\newcommand{\cT}{\mathcal{T}}
\newcommand{\cU}{\mathcal{U}}
\newcommand{\cB}{\mathcal{B}}
\newcommand{\Kb}{\overline{K_q}}

\begin{document}
\begin{sloppypar}

\title{The Maximum Second Zagreb Spectral Radius of Trees, Unicyclic Graphs and Bicyclic Graphs\thanks{This research was supported by the Higher Education Scientific Research Project of the Anhui Provincial Department of Education under Grants 2025AHGXZK31124, 2025AHGXZK31313 and SK2017A0586. ({\em Corresponding author: Genjie Wang})}}

\author{Juan Wang\thanks{Juan Wang is with the Faculty of Engineering, Anhui Sanlian University, Hefei 230601, Anhui, China. Email: wangjuan@mail.slu.edu.cn.}
\and Tony Newton\thanks{Tony Newton is with the Newton Astro Labs, London, United Kingdom. Email: tony.newton79@gmail.com.}
\and Genjie Wang\thanks{Genjie Wang is with the Faculty of Engineering, Anhui Sanlian University, Hefei 230601, Anhui, China. Email: gjwang2005@126.com.}
}
\date{}
\maketitle

\begin{abstract}
The second Zagreb matrix of a graph $G$ is $Z(G)=D(G)A(G)D(G)$, where $D(G)$ and $A(G)$ are the degree matrix and the adjacency matrix of $G$. For many degree-based weighted adjacency matrices, the trees, unicyclic graphs and bicyclic graphs with the largest spectral radius have a vertex adjacent to all others, and Mei et al. conjectured this for every such matrix on unicyclic graphs. Jiang et al. showed that $Z(G)$ is a counterexample. We prove that, for $q\in\{0,1,2\}$ and every $n\ge4$, the connected graph of order $n$ and cyclomatic number $q$ maximizing the spectral radius of $Z(G)$ is unique, and is obtained from two adjacent vertices with $q$ common neighbours by attaching $n-q-2$ pendent vertices to them as evenly as possible. The only exceptions are the star for $q=0$ and $n\le7$, and the triangle with one pendent vertex at each vertex for $q=1$ and $n=6$. This confirms a conjecture of Jiang et al.\ and all three parts of a conjecture of Shen and Shan for this matrix, with best possible thresholds on $n$. For $n\ge8$, the maximum value is $\lfloor (n+q)^2/4\rfloor+(n+q)/2+7q-2+O(1/n)$, whereas the star gives only $(n-1)^{3/2}$. A key tool is a determinant criterion that reduces spectral comparisons to sign conditions on explicit polynomials.

\vspace*{.3cm}
\end{abstract}

\section{Introduction}\label{sec:intro}
All graphs in this paper are finite, simple and connected. For a graph $G$, we write $d_v$ for the degree of a vertex $v$, and $D(G)$ and $A(G)$ for the diagonal degree matrix and the adjacency matrix of $G$. The first and second Zagreb indices go back to \cite{GutmanTrinajstic1972,GutmanRuscic1975}. The second Zagreb index and the corresponding second Zagreb matrix \cite{JafariRad2018} are
\[
M_2(G)=\sum_{uv\in E(G)}d_ud_v,\quad Z(G)=D(G)A(G)D(G).
\]
The $(u,v)$-entry of $Z(G)$ is $d_ud_v$ if $u$ and $v$ are adjacent and $0$ otherwise, so the sum of all entries of $Z(G)$ is $2M_2(G)$. We write $\rho(G)$ for the largest eigenvalue of $Z(G)$ and call it the second Zagreb spectral radius of $G$. Since $Z(G)$ is nonnegative and symmetric, $\rho(G)$ is also its spectral radius. Das, Mondal and Pal \cite{DasMondalPal2024} gave bounds for $\rho(G)$ and for the spread of $Z(G)$, and studied the energy of $Z(G)$. Relations between the Zagreb indices and the adjacency spectral radius or the energy of $G$ were studied, for example, in \cite{Abdo2014,DasGhalavand2025,Doslic2025}.

The matrix $Z(G)$ is a degree-based matrix, that is, its $(u,v)$-entry is $f(d_u,d_v)$ for adjacent $u,v$ and $0$ otherwise, where $f$ is a symmetric function \cite{Das2018,LiWang2021}. The second Zagreb matrix corresponds to $f(x,y)=xy$. Li and Wang \cite{LiWang2021} proposed a unified treatment of the spectral radii of such matrices for trees, and Li and Yang \cite{LiYang2024} gave unified bounds for their largest eigenvalue and their energy. Extremal unicyclic graphs and bicyclic graphs were then determined for several weights, for example for the first Zagreb matrix \cite{Das2024}, for the Sombor matrix \cite{Mei2023,Zhang2025}, for the atom-bond sum-connectivity matrix \cite{LinLiu2025}, and for weights with a property called $P^*$ \cite{Zheng2023,Ye2023}. For these weights, the maximizers contain a vertex of degree $n-1$. Mei, Fu, Miao and Gao \cite{Mei2023} conjectured that, for every degree-based matrix, the spectral radius of a unicyclic graph of order $n\ge3$ lies between the values for the cycle $C_n$ and for the graph obtained from a triangle by attaching $n-3$ pendent vertices to one of its vertices. Jiang, Zuo, Wang and Jahanbani \cite{Jiang2027} showed that this upper bound fails for the second Zagreb matrix. Based on exhaustive computations for $5\le n\le11$, they conjectured the extremal graphs described in Theorem~\ref{thm:main}(ii) below.

Shen and Shan \cite{Shen2025} studied positive symmetric weights $f$ that are increasing and convex in each variable, a class that contains $f(x,y)=xy$. By means of the Kelmans transformation \cite{Kelmans1981}, they obtained forbidden induced subgraphs for the graphs with the largest spectral radius among connected graphs of given order and cyclomatic number. They also introduced a property $P^{**}$, which $f(x,y)=xy$ satisfies, and conjectured \cite[Conjecture~6.7]{Shen2025} that for such $f$ and all sufficiently large $n$, the maximizers among trees, unicyclic graphs and bicyclic graphs have a common shape. In each of them, all vertices outside a small subgraph are pendent vertices attached to two adjacent vertices, and they are divided between these two vertices as evenly as possible. In this paper, we prove this conjecture for the second Zagreb matrix in all three classes, together with the exact range of orders and the exceptions for small orders. Part (ii) of our main result also confirms the conjecture of \cite{Jiang2027}.

To state the results, let $\cT_n$, $\cU_n$ and $\cB_n$ be the sets of trees, unicyclic graphs and bicyclic graphs of order $n$, that is, of connected graphs of order $n$ with cyclomatic number $0$, $1$ and $2$, respectively. For $q\ge0$, the book graph $K_2\vee\Kb$ consists of two adjacent vertices $u$ and $v$, which we call spine vertices, and $q$ further vertices adjacent to both $u$ and $v$. For $n\ge q+2$, let $H_{n,q}$ be the graph of order $n$ obtained from $K_2\vee\Kb$ by attaching $\lceil (n-q-2)/2\rceil$ pendent vertices to $u$ and $\lfloor (n-q-2)/2\rfloor$ pendent vertices to $v$. Thus, $H_{n,0}$ is the balanced double star, $H_{n,1}$ is obtained from a triangle by attaching $n-3$ pendent vertices to two of its vertices as evenly as possible, and $H_{n,2}$ is obtained from $K_4-e$ by attaching $n-4$ pendent vertices to its two vertices of degree $3$ as evenly as possible. For integers $a,b,c\ge0$, let $T(a,b,c)$ be the graph obtained from a triangle by attaching $a$, $b$ and $c$ pendent vertices to its three vertices, respectively. It has order $a+b+c+3$, and up to isomorphism it does not depend on the order of $a,b,c$. Let $S_n$ denote the star of order $n$. Some of these graphs are drawn in Figure~\ref{fig:graphs}.

\begin{theorem}\label{thm:main}
	The following statements hold.
	\begin{itemize}
		\item[(i)] Let $n\ge4$ and $T\in\cT_n$. If $n\le7$, then $\rho(T)\le\rho(S_n)$, with equality if and only if $T\cong S_n$. If $n\ge8$, then $\rho(T)\le\rho(H_{n,0})$, with equality if and only if $T\cong H_{n,0}$.
		\item[(ii)] Let $n\ge3$ and $G\in\cU_n$. Then $\rho(G)\le\rho(G_n^*)$, with equality if and only if $G\cong G_n^*$, where $G_6^*=T(1,1,1)$ and $G_n^*=H_{n,1}$ for $n\ne6$.
		\item[(iii)] Let $n\ge4$ and $G\in\cB_n$. Then $\rho(G)\le\rho(H_{n,2})$, with equality if and only if $G\cong H_{n,2}$.
	\end{itemize}
\end{theorem}

Theorem~\ref{thm:main}(ii) for $n\ge5$ is Conjecture~1 of \cite{Jiang2027}. For $f(x,y)=xy$, parts (i), (ii) and (iii) of Theorem~\ref{thm:main} prove parts (1), (2) and (3) of \cite[Conjecture~6.7]{Shen2025} for $n\ge8$, $n\ge7$ and $n\ge4$, respectively. The exceptions in (i) and (ii) show that the thresholds $8$ and $7$ cannot be lowered. The corresponding minimum problems are known. For $f(x,y)=xy$, the path $P_n$ is the unique tree with the smallest $\rho$ for $n\ge9$ \cite[Theorem~1.6]{LiWang2021}, the cycle $C_n$ is the unique unicyclic graph with the smallest $\rho$, as follows from \cite[Theorem~5.1]{Shen2025}, and the bicyclic graphs with the smallest $\rho$ were determined for $n\ge8$ in \cite[Theorem~5.6]{Shen2025}.

Our second result gives the maximum values. The quartic polynomials $F^{(q)}_{s,p}$ used in part (ii) are defined in \eqref{eq:F} below.

\begin{theorem}\label{thm:values}
	Let $q\in\{0,1,2\}$, $n\ge q+2$ and $s=n+q$.
	\begin{itemize}
		\item[(i)] If $s=2t$ is even, then
		\[
		\rho(H_{n,q})=\frac t2\left(t+\sqrt{t^2+4t+28q-4}\right).
		\]
		\item[(ii)] If $s$ is odd, then $\rho(H_{n,q})$ is the largest real root of $F^{(q)}_{s,(s^2-1)/4}(\lambda)$.
		\item[(iii)] If $s\ge7$, then
		\[
		\left|\rho(H_{n,q})-\left\lfloor\frac{s^2}{4}\right\rfloor-\frac s2-7q+2\right|<\frac{64}{s}.
		\]
	\end{itemize}
	Moreover, $\rho(S_n)=(n-1)^{3/2}$ and $\rho(T(1,1,1))=9+3\sqrt{10}$.
\end{theorem}

By Theorems~\ref{thm:main} and~\ref{thm:values}, the maximum of $\rho$ over connected graphs of order $n\ge8$ with cyclomatic number $q\in\{0,1,2\}$ is $\lfloor (n+q)^2/4\rfloor+(n+q)/2+7q-2+O(1/n)$. In comparison, $\rho(S_n)=(n-1)^{3/2}$. Table~\ref{tab:max} lists the maximum values for small orders.

\begin{table}[ht]
	\centering
	\caption{Maximum values of $\rho$ for $4\le n\le20$ (rounded).}\label{tab:max}
	\begin{tabular}{rrrr@{\qquad}rrrr}
		\toprule
		$n$ & $\cT_n$ & $\cU_n$ & $\cB_n$ & $n$ & $\cT_n$ & $\cU_n$ & $\cB_n$\\
		\midrule
		4 & 5.196152 & 11.040986 & 17.316006 & 13 & 47.054209 & 59.674565 & 72.086949 \\
		5 & 8.000000 & 14.562306 & 21.550649 & 14 & 54.404013 & 67.366717 & 80.662100 \\
		6 & 11.180340 & 18.486833 & 26.330303 & 15 & 61.997849 & 75.817805 & 89.447401 \\
		7 & 14.696938 & 22.966630 & 31.419138 & 16 & 70.366652 & 84.484389 & 99.000000 \\
		8 & 18.583005 & 27.796785 & 37.122145 & 17 & 78.952061 & 93.934539 & 108.746322 \\
		9 & 23.218445 & 33.266560 & 43.102260 & 18 & 88.335656 & 103.580992 & 119.282032 \\
		10 & 28.507811 & 39.036133 & 49.749016 & 19 & 97.914088 & 114.031242 & 129.997102 \\
		11 & 34.125405 & 45.495454 & 56.646832 & 20 & 108.309519 & 124.661530 & 141.520059 \\
		12 & 40.449944 & 52.220809 & 64.252358 & & & & \\
		\bottomrule
	\end{tabular}
\end{table}

The proof of Theorem~\ref{thm:main} starts from the reductions of Li and Wang \cite{LiWang2021} and of Shen and Shan \cite{Shen2025}. These reductions show that every maximizer is obtained from $K_2$, $C_3$, $K_4-e$ or $K_{2,3}$ by attaching pendent vertices in a restricted way. All four graphs are complete multipartite graphs. For graphs of this type, we prove in Lemma~\ref{lem:det} that the sign of an explicit determinant decides whether $\rho$ is smaller than a given number. The comparisons are then made by a balancing lemma for the graphs with two adjacent hubs, by concavity, and by polynomial inequalities. Several of these inequalities are proved by direct estimates, the remaining ones are certified by expansions with positive integer coefficients, and the small orders are settled by exact rational computations. All these computations use exact integer and rational arithmetic.

The paper is organized as follows. Section~\ref{sec:prelim} recalls the reductions and proves the determinant criterion. Section~\ref{sec:book} treats the graphs with two adjacent hubs, proves Theorem~\ref{thm:values}, and ends with the proof of Theorem~\ref{thm:main}(i) for trees. Sections~\ref{sec:uni} and~\ref{sec:bi} prove parts (ii) and (iii) of Theorem~\ref{thm:main}. Section~\ref{sec:concl} contains concluding remarks and a conjecture for larger cyclomatic numbers.

\section{Reductions and a determinant criterion}\label{sec:prelim}
\subsection{Reductions}
For a connected graph $G$ with $n$ vertices and $m$ edges, $m-n+1$ is its cyclomatic number. A vertex of degree one is pendent. If $G$ contains a cycle, the base $\widehat G$ of $G$ is the graph obtained from $G$ by repeatedly deleting pendent vertices. The function $f(x,y)=xy$ is positive and symmetric for $x,y>0$, and it is increasing and convex in each variable. Hence, the following results of \cite{LiWang2021} and \cite{Shen2025} apply to the second Zagreb matrix.

\begin{lemma}[{\cite[Theorem~1.5]{LiWang2021}}]\label{lem:trees}
	Every graph with the largest $\rho$ in $\cT_n$ is the star $S_n$ or a double star.
\end{lemma}

\begin{lemma}[{\cite[Theorem~6.2 and Corollaries~6.3--6.5]{Shen2025}}]\label{lem:shen}
	Let $q\in\{1,2\}$, and let $G$ be a graph with the largest $\rho$ among connected graphs of order $n$ with cyclomatic number $q$. Then the following hold.
	\begin{itemize}
		\item[(a)] $G$ has no induced subgraph isomorphic to the path $P_5$.
		\item[(b)] Every vertex of $G$ outside $\widehat G$ is pendent.
		\item[(c)] If $q=1$, then $\widehat G=C_3$. If $q=2$, then $\widehat G$ is isomorphic to $K_4-e$ or $K_{2,3}$.
	\end{itemize}
\end{lemma}

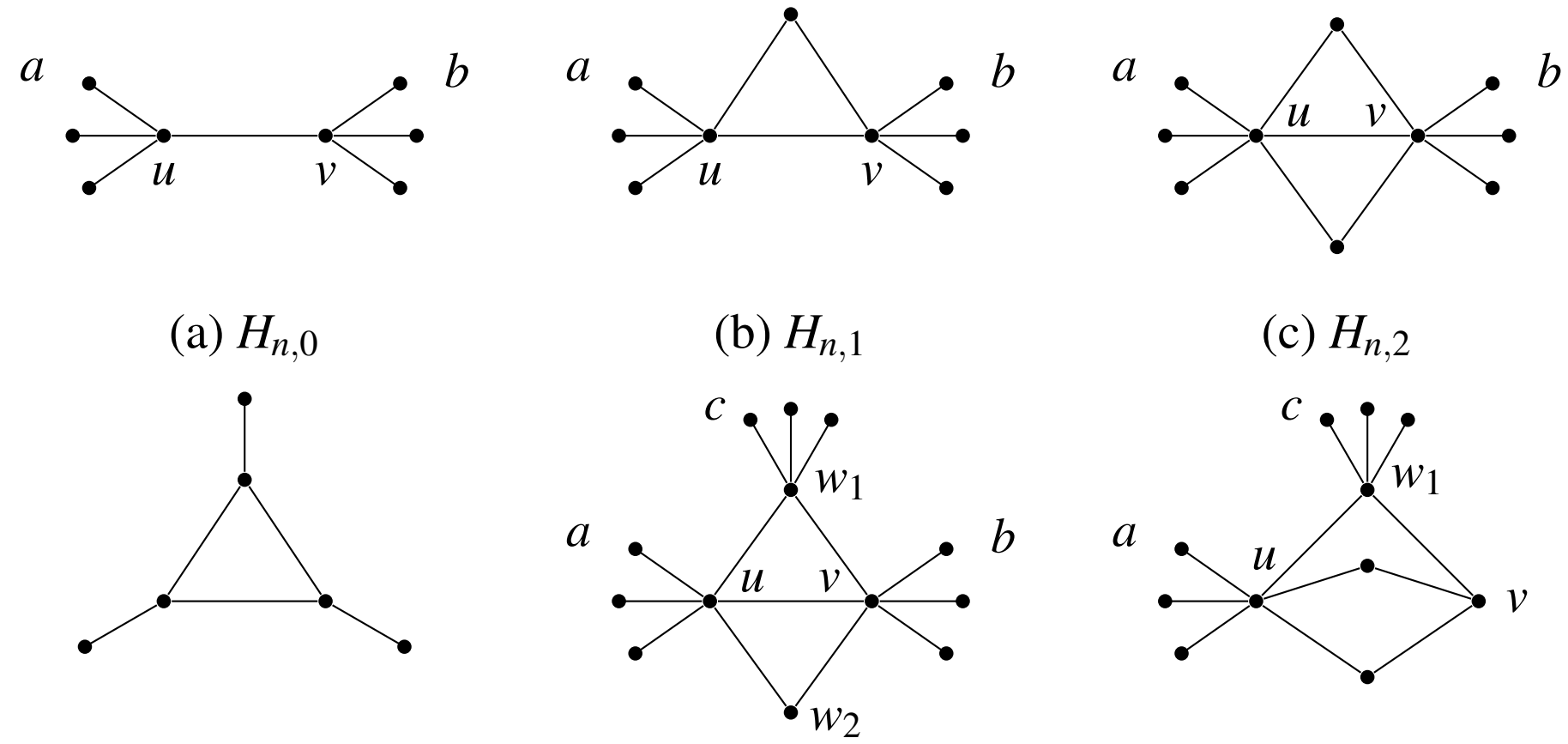
\begin{figure}[ht]
	\centering
	\begin{tikzpicture}[v/.style={circle,fill=black,inner sep=1.1pt},scale=0.78,every node/.style={font=\small}]
		\begin{scope}[shift={(0,0)}]
			\node[v,label={below:$u$}] (u) at (0,0) {};
			\node[v,label={below:$v$}] (v) at (1.6,0) {};
			\draw (u)--(v);
			\foreach \t in {145,180,215} \draw (u)--++(\t:0.9) node[v]{};
			\foreach \t in {35,0,-35} \draw (v)--++(\t:0.9) node[v]{};
			\node at (-1.3,0.65) {$a$}; \node at (2.9,0.65) {$b$};
			\node at (0.8,-2.0) {(a) $H_{n,0}$};
		\end{scope}
		\begin{scope}[shift={(5.4,0)}]
			\node[v,label={below:$u$}] (u) at (0,0) {};
			\node[v,label={below:$v$}] (v) at (1.6,0) {};
			\node[v] (w) at (0.8,1.2) {};
			\draw (u)--(v)--(w)--(u);
			\foreach \t in {145,180,215} \draw (u)--++(\t:0.9) node[v]{};
			\foreach \t in {35,0,-35} \draw (v)--++(\t:0.9) node[v]{};
			\node at (-1.3,0.65) {$a$}; \node at (2.9,0.65) {$b$};
			\node at (0.8,-2.0) {(b) $H_{n,1}$};
		\end{scope}
		\begin{scope}[shift={(10.8,0)}]
			\node[v] (u) at (0,0) {}; \node at (0.42,0.2) {$u$};
			\node[v] (v) at (1.6,0) {}; \node at (1.18,0.2) {$v$};
			\node[v] (w1) at (0.8,1.1) {};
			\node[v] (w2) at (0.8,-1.1) {};
			\draw (u)--(v); \draw (u)--(w1)--(v); \draw (u)--(w2)--(v);
			\foreach \t in {145,180,215} \draw (u)--++(\t:0.9) node[v]{};
			\foreach \t in {35,0,-35} \draw (v)--++(\t:0.9) node[v]{};
			\node at (-1.3,0.65) {$a$}; \node at (2.9,0.65) {$b$};
			\node at (0.8,-2.0) {(c) $H_{n,2}$};
		\end{scope}
		\begin{scope}[shift={(0,-4.6)}]
			\node[v] (u) at (0,0) {};
			\node[v] (v) at (1.6,0) {};
			\node[v] (w) at (0.8,1.2) {};
			\draw (u)--(v)--(w)--(u);
			\draw (u)--++(210:0.9) node[v]{};
			\draw (v)--++(-30:0.9) node[v]{};
			\draw (w)--++(90:0.8) node[v]{};
			\node at (0.8,-2.0) {(d) $T(1,1,1)$};
		\end{scope}
		\begin{scope}[shift={(5.4,-4.6)}]
			\node[v] (u) at (0,0) {}; \node at (0.42,0.2) {$u$};
			\node[v] (v) at (1.6,0) {}; \node at (1.18,0.2) {$v$};
			\node[v] (w1) at (0.8,1.1) {}; \node at (1.3,1.15) {$w_1$};
			\node[v] (w2) at (0.8,-1.1) {}; \node at (1.25,-1.2) {$w_2$};
			\draw (u)--(v); \draw (u)--(w1)--(v); \draw (u)--(w2)--(v);
			\foreach \t in {145,180,215} \draw (u)--++(\t:0.9) node[v]{};
			\foreach \t in {35,0,-35} \draw (v)--++(\t:0.9) node[v]{};
			\foreach \t in {60,90,120} \draw (w1)--++(\t:0.8) node[v]{};
			\node at (-1.3,0.65) {$a$}; \node at (2.9,0.65) {$b$}; \node at (0.05,1.9) {$c$};
			\node at (0.8,-2.0) {(e) $D(a,b,c)$};
		\end{scope}
		\begin{scope}[shift={(10.8,-4.6)}]
			\node[v] (u) at (0,0) {}; \node at (0.08,0.42) {$u$};
			\node[v,label={right:$v$}] (v) at (2.2,0) {};
			\node[v] (w1) at (1.1,1.1) {}; \node at (1.6,1.2) {$w_1$};
			\node[v] (w2) at (1.1,0.35) {};
			\node[v] (w3) at (1.1,-0.75) {};
			\foreach \w in {w1,w2,w3} {\draw (u)--(\w)--(v);}
			\foreach \t in {145,180,215} \draw (u)--++(\t:0.9) node[v]{};
			\foreach \t in {60,90,120} \draw (w1)--++(\t:0.8) node[v]{};
			\node at (-1.3,0.65) {$a$}; \node at (0.35,1.9) {$c$};
			\node at (1.1,-2.0) {(f) $K(a,c)$};
		\end{scope}
	\end{tikzpicture}
	\caption{The graphs $H_{n,q}$ for $q=0,1,2$, the graph $T(1,1,1)$, and the graphs $D(a,b,c)$ and $K(a,c)$. A fan of three edges labelled $a$, $b$ or $c$ stands for that number of pendent vertices. In (a)--(c), $a=\lceil(n-q-2)/2\rceil$ and $b=\lfloor(n-q-2)/2\rfloor$.}\label{fig:graphs}
\end{figure}

In \cite{Shen2025}, the graphs $K_4-e$ and $K_{2,3}$ are denoted by $\Theta(1,2,2)$ and $\Theta(2,2,2)$. The results of \cite[Section~6]{Shen2025} are derived from forbidden induced subgraphs obtained with the Kelmans transformation \cite{Kelmans1981}, whose effect on the spectral radius was studied by Csikv\'ari \cite{Csikvari2009} for the adjacency matrix and by Li and Wang \cite{LiWang2021} for degree-based weighted adjacency matrices.
For integers $a,b,c\ge0$, let $D(a,b,c)$ be the graph obtained from $K_4-e$ by attaching $a$ and $b$ pendent vertices to the two vertices of degree $3$ and $c$ pendent vertices to one vertex of degree $2$. Let $K(a,c)$ be the graph obtained from $K_{2,3}$ by attaching $a$ pendent vertices to one vertex of degree $3$ and $c$ pendent vertices to one vertex of degree $2$. These graphs have orders $a+b+c+4$ and $a+c+5$, respectively, and they are drawn in Figure~\ref{fig:graphs}.

\begin{corollary}\label{cor:families}
	Every graph with the largest $\rho$ in $\cU_n$ is isomorphic to $T(a,b,c)$ for some $a,b,c\ge0$. Every graph with the largest $\rho$ in $\cB_n$ is isomorphic to $D(a,b,c)$ or to $K(a,c)$ for some $a,b,c\ge0$.
\end{corollary}

\begin{proof}
	The first statement follows from Lemma~\ref{lem:shen}(b) and (c). Let $G$ be a graph with the largest $\rho$ in $\cB_n$. By Lemma~\ref{lem:shen}(b), $G$ is obtained from $\widehat G$ by attaching pendent vertices. Suppose first that $\widehat G=K_4-e$, where $u,v$ have degree $3$ and $w_1,w_2$ have degree $2$ in $\widehat G$. If $w_1$ and $w_2$ carry pendent vertices $\ell_1$ and $\ell_2$, respectively, then $\ell_1w_1uw_2\ell_2$ is an induced path, since $w_1$ and $w_2$ are not adjacent. This contradicts Lemma~\ref{lem:shen}(a), so $G$ is isomorphic to some $D(a,b,c)$. Suppose next that $\widehat G=K_{2,3}$. If two distinct vertices $y_1,y_2$ of the same part carry pendent vertices $\ell_1$ and $\ell_2$, respectively, choose a vertex $y_3$ of the other part. Then $\ell_1y_1y_3y_2\ell_2$ is an induced path, again a contradiction. Hence, at most one vertex of each part carries pendent vertices, and $G$ is isomorphic to some $K(a,c)$.
\end{proof}

\subsection{A determinant criterion}
The following lemma replaces the computation of $\rho$ by the sign of a determinant whose size is the order of the base. For a vector $\bd$ and a set $V_t$ of indices, $\bd_t$ denotes the vector obtained from $\bd$ by replacing the entries outside $V_t$ with $0$.

\begin{lemma}\label{lem:det}
	Let $\Gamma$ be a connected complete multipartite graph with vertex set $\{1,\dots,k\}$, where $k\ge2$, and with parts $V_1,\dots,V_r$. Let $G$ be obtained from $\Gamma$ by attaching $p_i\ge0$ pendent vertices to each vertex $i$ of $\Gamma$, and put $d_i=d_G(i)$, $D=\diag(d_1,\dots,d_k)$ and $\bd=(d_1,\dots,d_k)^{\T}$. For $\lambda>0$, define
	\begin{align*}
		M_G(\lambda)&=\diag\big(\lambda^2-p_1d_1^2,\dots,\lambda^2-p_kd_k^2\big)-\lambda DA(\Gamma)D,\\
		N_G(\lambda)&=\diag\big(\lambda^2-p_1d_1^2,\dots,\lambda^2-p_kd_k^2\big)+\lambda\sum_{t=1}^{r}\bd_t\bd_t^{\T}.
	\end{align*}
	If $N_G(\lambda)$ is positive definite, then
	\begin{equation}\label{eq:mdl}
		\det M_G(\lambda)=\det N_G(\lambda)\big(1-\lambda\bd^{\T}N_G(\lambda)^{-1}\bd\big),
	\end{equation}
	and $\det M_G(\lambda)$ is positive, zero or negative according as $\rho(G)<\lambda$, $\rho(G)=\lambda$ or $\rho(G)>\lambda$. Moreover, $N_G(\lambda)$ is positive definite if $\lambda>p_i$ for every vertex $i$ that forms a part of size one, and $\lambda^2>p_id_i^2$ for every vertex $i$ in a part of size at least two.
\end{lemma}

\begin{proof}
	Order the vertices of $G$ so that the vertices of $\Gamma$ come first. Each of the $p_i$ pendent vertices at $i$ is joined to $i$ by an edge of weight $d_i$. Hence,
	\[
	\lambda I-Z(G)=\begin{pmatrix}\lambda I_k-DA(\Gamma)D&-E\\-E^{\T}&\lambda I\end{pmatrix},
	\]
	where each column of $E$ has exactly one nonzero entry, namely $d_i$ in row $i$ for a pendent vertex at $i$. Thus, $EE^{\T}=\diag(p_1d_1^2,\dots,p_kd_k^2)$, and the Schur complement of the block $\lambda I$ is $\lambda I_k-DA(\Gamma)D-\lambda^{-1}EE^{\T}=\lambda^{-1}M_G(\lambda)$. Since $\lambda I$ is positive definite, Haynsworth's inertia additivity formula \cite{Haynsworth1968} shows that $\lambda I-Z(G)$ and $M_G(\lambda)$ have the same number of negative eigenvalues and the same number of zero eigenvalues. If $G=\Gamma$, this holds trivially, because $M_G(\lambda)=\lambda(\lambda I-Z(G))$.

	Two vertices of $\Gamma$ are adjacent if and only if they lie in different parts. Hence, $A(\Gamma)=J-\sum_t\mathbf 1_t\mathbf 1_t^{\T}$, where $J$ is the all-ones matrix and $\mathbf 1_t$ is the indicator vector of $V_t$. Therefore, $DA(\Gamma)D=\bd\bd^{\T}-\sum_t\bd_t\bd_t^{\T}$, that is,
	\[
	M_G(\lambda)=N_G(\lambda)-\lambda\bd\bd^{\T}.
	\]
	Now assume that $N_G(\lambda)$ is positive definite. Then \eqref{eq:mdl} is the matrix determinant lemma. Let $\mu_1\le\mu_2\le\dots\le\mu_k$ be the eigenvalues of $M_G(\lambda)$. By the interlacing theorem for rank-one perturbations \cite[Section~4.3]{HornJohnson2013}, $\mu_2$ is at least the smallest eigenvalue of $N_G(\lambda)$, so $\mu_2>0$. Hence, $M_G(\lambda)$ has at most one nonpositive eigenvalue, and $\det M_G(\lambda)$ has the sign of $\mu_1$. If $\mu_1>0$, then $\lambda I-Z(G)$ is positive definite, that is, $\rho(G)<\lambda$. If $\mu_1=0$, then $\lambda I-Z(G)$ is positive semidefinite and singular, so $\lambda=\rho(G)$. If $\mu_1<0$, then $\lambda I-Z(G)$ has a negative eigenvalue, so $\rho(G)>\lambda$.

	Finally, $N_G(\lambda)$ is block diagonal with respect to the parts. The block of a part $\{i\}$ is the number $\lambda^2+(\lambda-p_i)d_i^2$. The block of a part $V_t$ with $|V_t|\ge2$ is the sum of the diagonal matrix with entries $\lambda^2-p_id_i^2$, $i\in V_t$, and a positive semidefinite matrix. This proves the last statement.
\end{proof}

For $\Gamma=K_3$, all parts have size one and $N_G(\lambda)=\diag(h_1,h_2,h_3)$ with $h_i=\lambda^2+(\lambda-p_i)d_i^2$. Then \eqref{eq:mdl} reads $\det M_G(\lambda)=h_1h_2h_3\big(1-\sum_{i=1}^3\lambda d_i^2/h_i\big)$, which is the usual secular equation obtained by eliminating the pendent vertices of a triangle.
We also use the Perron--Frobenius theorem in the following form. If $G$ is connected, then $Z(G)$ is irreducible, and every proper principal submatrix of $Z(G)$ has a smaller spectral radius than $Z(G)$ \cite[Chapter~8]{HornJohnson2013}. In particular,
\begin{equation}\label{eq:pf}
	\rho(G)>d_ud_v\quad\text{for every edge $uv$ of a connected graph $G$ of order at least $3$.}
\end{equation}

\section{Graphs with two adjacent hubs}\label{sec:book}
In this section, $q\in\{0,1,2\}$. For integers $x,y\ge q+1$, let $B_q[x,y]$ be the graph obtained from $K_2\vee\Kb$, with spine vertices $u,v$ and further vertices $w_1,\dots,w_q$, by attaching $x-q-1$ pendent vertices to $u$ and $y-q-1$ pendent vertices to $v$. Then $d_u=x$, $d_v=y$, and the order of $B_q[x,y]$ is $x+y-q$. Throughout, $s=x+y$ and $p=xy$. Thus, $B_0[x,y]$ is a double star for $y\ge2$, $B_0[n-1,1]=S_n$, $B_1[x,y]=T(x-2,y-2,0)$, $B_2[x,y]=D(x-3,y-3,0)$, and $H_{n,q}=B_q[\lceil s/2\rceil,\lfloor s/2\rfloor]$ with $s=n+q$. For a real number $\kappa$ we put
\begin{equation}\label{eq:lamk}
	\lambda_\kappa(s)=\frac{s^2+2s+\kappa}{4}.
\end{equation}

\begin{lemma}\label{lem:book}
	Let $x\ge y\ge q+1$ and $\lambda>x-q-1$. Then $\det M_{B_q[x,y]}(\lambda)=\lambda^{2q}F^{(q)}_{s,p}(\lambda)$, where
	\begin{equation}\label{eq:F}
		\begin{aligned}
			F^{(q)}_{s,p}(\lambda)={}&\lambda^4-\big[p^2-(3s+6q-2)p+s^3+(3q-1)s^2\big]\lambda^2-8qp^2\lambda\\
			&+p^2\big[p+(3q-1)s+(3q-1)^2-16q^2\big].
		\end{aligned}
	\end{equation}
	In particular, $F^{(q)}_{s,p}(\lambda)$ has the sign of $\lambda-\rho(B_q[x,y])$. Moreover,
	\begin{equation}\label{eq:Ffactor}
		F^{(q)}_{2t,t^2}(\lambda)=\big[\lambda^2-t^2\lambda-(t+7q-1)t^2\big]\big[\lambda^2+t^2\lambda-(t-q-1)t^2\big].
	\end{equation}
\end{lemma}

\begin{proof}
	The graph $K_2\vee\Kb$ is complete multipartite with parts $\{u\}$, $\{v\}$ and, if $q\ge1$, $\{w_1,\dots,w_q\}$. The vertices $w_i$ carry no pendent vertices, and $\lambda>x-q-1\ge y-q-1$. Hence, Lemma~\ref{lem:det} applies. In $M_{B_q[x,y]}(\lambda)$, the diagonal entries at $u$, $v$ and $w_i$ are $\lambda^2-(x-q-1)x^2$, $\lambda^2-(y-q-1)y^2$ and $\lambda^2$. The off-diagonal entries are $-\lambda xy$ at $uv$, $-2\lambda x$ at $uw_i$, $-2\lambda y$ at $vw_i$, and $0$ between distinct vertices $w_i,w_j$. Taking the Schur complement of the block $\lambda^2I_q$ indexed by $w_1,\dots,w_q$, we obtain
	\[
	\det M_{B_q[x,y]}(\lambda)=\lambda^{2q}\det\begin{pmatrix}\lambda^2-(x+3q-1)x^2&-(\lambda+4q)xy\\-(\lambda+4q)xy&\lambda^2-(y+3q-1)y^2\end{pmatrix}.
	\]
	Expanding the determinant with $x^2+y^2=s^2-2p$ and $x^3+y^3=s^3-3sp$ gives \eqref{eq:F}. For $x=y=t$, the determinant is the difference of the squares of $\lambda^2-(t+3q-1)t^2$ and $(\lambda+4q)t^2$, which gives \eqref{eq:Ffactor}.
\end{proof}

For later reference, the three polynomials in \eqref{eq:F} are
\begin{align*}
	F^{(0)}_{s,p}(\lambda)&=\lambda^4-\big[p^2-(3s-2)p+s^3-s^2\big]\lambda^2+p^2(p-s+1),\\
	F^{(1)}_{s,p}(\lambda)&=\lambda^4-\big[p^2-(3s+4)p+s^3+2s^2\big]\lambda^2-8p^2\lambda+p^2(p+2s-12),\\
	F^{(2)}_{s,p}(\lambda)&=\lambda^4-\big[p^2-(3s+10)p+s^3+5s^2\big]\lambda^2-16p^2\lambda+p^2(p+5s-39).
\end{align*}

The next lemma shows that, for fixed $s$, the spectral radius of $B_q[x,y]$ increases as $x$ and $y$ become closer.

\begin{lemma}\label{lem:balance}
	Let $s_0=9$, $s_1=12$, $s_2=10$, $m_0=m_1=2$ and $m_2=3$. Let $s\ge s_q$, and let $x\ge y\ge m_q$ and $x'\ge y'\ge m_q$ be integers with $x+y=x'+y'=s$ and $xy<x'y'$. Then $\rho(B_q[x,y])<\rho(B_q[x',y'])$. In particular, $B_q[\lceil s/2\rceil,\lfloor s/2\rfloor]$ is the unique graph with the largest $\rho$ among the graphs $B_q[x,y]$ with $x+y=s$ and $x\ge y\ge m_q$.
\end{lemma}

\begin{proof}
	Put $\lambda'=\rho(B_q[x',y'])$, $p=xy$ and $p'=x'y'$. The edge $uv$ of $B_q[x',y']$ has weight $p'$, so $\lambda'>p'$ by \eqref{eq:pf}. Since $p'\ge m_q(s-m_q)\ge2s-4>s$, we have $\lambda'>x-q-1$ and $\lambda'>x'-q-1$. By Lemma~\ref{lem:book}, $F^{(q)}_{s,p'}(\lambda')=0$, and it suffices to show that $F^{(q)}_{s,p}(\lambda')>0$. We regard $F^{(q)}_{s,\tilde p}(\lambda')$ as a polynomial in a real variable $\tilde p\in[m_q(s-m_q),p']$, which contains $p$. Its derivative is
	\begin{equation}\label{eq:dF}
		\frac{\partial F^{(q)}_{s,\tilde p}(\lambda')}{\partial \tilde p}=-(2\tilde p-3s-6q+2)\lambda'^2-16q\tilde p\lambda'+3\tilde p^2+2\big[(3q-1)s+(3q-1)^2-16q^2\big]\tilde p.
	\end{equation}
	For $\tilde p\ge m_q(s-m_q)$ and $s\ge s_q$, we have $2\tilde p\ge3s+6q-2$, since for $q=0,1,2$ this follows from $4s-8\ge3s-2$, $4s-8\ge3s+4$ and $6s-18\ge3s+10$, respectively. Hence, the first term in \eqref{eq:dF} is nonpositive. Recall also that $\lambda'>p'\ge\tilde p$.

	If $q\ge1$, then dropping the first term and using $\lambda'>\tilde p$ in the second term gives
	\[
	\frac{\partial F^{(q)}_{s,\tilde p}(\lambda')}{\partial \tilde p}<\tilde p\big[(3-16q)\tilde p+2(3q-1)s+2(3q-1)^2-32q^2\big].
	\]
	For $q=1$, the bracket is $-13\tilde p+4s-24<0$ because $\tilde p\ge2s-4$. For $q=2$, it is $-29\tilde p+10s-78<0$ because $\tilde p\ge3s-9$. If $q=0$, then using $\lambda'^2>\tilde p^2$ in the first term gives
	\[
	\frac{\partial F^{(0)}_{s,\tilde p}(\lambda')}{\partial \tilde p}\le-\tilde p\big[2\tilde p^2-(3s+1)\tilde p+2(s-1)\big]<0,
	\]
	because $2\tilde p\ge4s-8\ge3s+1$ for $s\ge9$. Thus, the derivative is negative in all cases, so $F^{(q)}_{s,p}(\lambda')>F^{(q)}_{s,p'}(\lambda')=0$. By Lemma~\ref{lem:book}, $\rho(B_q[x,y])<\lambda'$.
\end{proof}

We next give lower bounds for $\rho(H_{n,1})$ and $\rho(H_{n,2})$ by numbers of the form \eqref{eq:lamk}. They will separate these graphs from all other candidates.

\begin{lemma}\label{lem:lower}
	For every $n\ge8$, we have
	\[
	\rho(H_{n,1})>\lambda_{12}(n+1)\quad\text{and}\quad \rho(H_{n,2})>\lambda_{27}(n+2).
	\]
\end{lemma}

\begin{proof}
	Let $(q,\kappa)\in\{(1,12),(2,27)\}$, $s=n+q$ and $\lambda=\lambda_\kappa(s)$. Since $\lambda>s^2/4>s>\lceil s/2\rceil$, Lemma~\ref{lem:book} applies to $H_{n,q}=B_q[\lceil s/2\rceil,\lfloor s/2\rfloor]$ at $\lambda$, and it suffices to show that $F^{(q)}_{s,\lfloor s^2/4\rfloor}(\lambda)<0$.
	Suppose first that $s=2t$ is even. Since $\lambda\ge t$, the second factor in \eqref{eq:Ffactor} is positive, and with $\lambda=t^2+t+\kappa/4$ the first factor equals
	\[
	\Big(\frac\kappa4+2-7q\Big)t^2+\frac\kappa2t+\frac{\kappa^2}{16}=
	\begin{cases}
		-2t^2+6t+9=-\big(2t(t-3)-9\big),&(q,\kappa)=(1,12),\\[2pt]
		-\frac{1}{16}\big(12t(7t-18)-729\big),&(q,\kappa)=(2,27).
	\end{cases}
	\]
	Both values are negative for $t\ge5$, since $2t(t-3)\ge20$ and $12t(7t-18)\ge1020$. For $q=1$ the condition $t\ge5$ means $n\ge9$, and for $q=2$ it means $n\ge8$.

	Suppose next that $s$ is odd, so that $\lfloor s^2/4\rfloor=(s^2-1)/4$. Direct substitution into \eqref{eq:F} gives
	\begin{align*}
		256F^{(1)}_{s,(s^2-1)/4}(\lambda_{12}(s))&=-14s^6+60s^5+139s^4+2252s^3+8432s^2+11248s+17708,\\
		256F^{(2)}_{s,(s^2-1)/4}(\lambda_{27}(s))&=-40s^6+128s^5-4s^4+10576s^3+69504s^2+144240s+499196.
	\end{align*}
	For $s\ge9$, the first polynomial is negative, because $14s^6-60s^5=2s^5(7s-30)\ge66s^5$, while $139s^4+2252s^3+8432s^2+11248s+17708<57s^5$. For $s\ge11$, the second polynomial is negative, because $40s^6-128s^5+4s^4>8s^5(5s-16)\ge312s^5$, while $10576s^3+69504s^2+144240s+499196<153s^5$. For $q=1$ the condition $s\ge9$ means $n\ge8$, and for $q=2$ the condition $s\ge11$ means $n\ge9$. Combining the two parities gives the lemma.
\end{proof}

\begin{proof}[Proof of Theorem~\ref{thm:values}]
	For $H_{2,0}=K_2$, we have $\rho=1$, which agrees with (i). Let $n\ge3$, let $x=\lceil s/2\rceil$ and $y=\lfloor s/2\rfloor$, so that $H_{n,q}=B_q[x,y]$, and recall that $\rho(H_{n,q})>xy\ge x$ by \eqref{eq:pf}. By Lemma~\ref{lem:book}, $F^{(q)}_{s,xy}(\lambda)>0$ for $\lambda>\rho(H_{n,q})$ and $F^{(q)}_{s,xy}(\rho(H_{n,q}))=0$. This proves (ii).

	For (i), let $s=2t$. For $\lambda\ge t$, the second factor in \eqref{eq:Ffactor} is positive. The first factor is negative at $\lambda=t$, and its larger root is
	$\lambda_+=\frac t2\left(t+\sqrt{t^2+4t+28q-4}\right)>t$.
	Hence, $F^{(q)}_{2t,t^2}$ is negative on $(t,\lambda_+)$ and positive on $(\lambda_+,\infty)$. Since $\rho(H_{n,q})>t^2\ge t$, we obtain $\rho(H_{n,q})=\lambda_+$.

	For (iii), let $s\ge7$, let $\varepsilon\in\{0,1\}$ be the parity of $s$, put $\beta=(s^2-\varepsilon)/4+s/2+7q-2$, and let $\lambda_\pm=\beta\pm64/s$. Since $x\le(s+1)/2$, we have $\lambda_--x\ge(s^2-11)/4-64/s>0$, because $s(s^2-11)\ge7\cdot38>256$. Thus, $\lambda_->x>x-q-1$, and Lemma~\ref{lem:book} applies at $\lambda_\pm$. Since $4s\lambda_\pm=s(s^2+2s+28q-8-\varepsilon)\pm256$, the numbers $256s^4F^{(q)}_{s,(s^2-\varepsilon)/4}(\lambda_\pm)$ are polynomials of degree $9$ in $s$, one for each choice of $q$, $\varepsilon$ and the sign. Substituting $s=3+\sigma$ into the six polynomials for $\lambda_+$, and $s=7+\sigma$ into the six polynomials for $\lambda_-$, yields polynomials in $\sigma$ whose coefficients are all positive for $\lambda_+$ and all negative for $\lambda_-$. Therefore, $F^{(q)}_{s,(s^2-\varepsilon)/4}(\lambda_+)>0>F^{(q)}_{s,(s^2-\varepsilon)/4}(\lambda_-)$, which gives $\lambda_-<\rho(H_{n,q})<\lambda_+$.

	Finally, $Z(S_n)$ has the eigenvalue $(n-1)^{3/2}$ with a positive eigenvector, so $\rho(S_n)=(n-1)^{3/2}$. In $T(1,1,1)$, the cycle vertices have degree $3$, and a positive eigenvector is constant on the cycle vertices and on the pendent vertices. With values $1$ and $\eta$ on these two classes, the eigenvalue equations are $\rho=18+3\eta$ and $\rho\eta=3$, so $\rho^2-18\rho-9=0$ and $\rho(T(1,1,1))=9+3\sqrt{10}$.
\end{proof}

We close this section with the proof of Theorem~\ref{thm:main}(i).

\begin{proof}[Proof of Theorem~\ref{thm:main}(i)]
	Let $T$ be a graph with the largest $\rho$ in $\cT_n$. By Lemma~\ref{lem:trees}, $T=B_0[x,y]$ for some integers $x\ge y\ge1$ with $x+y=n$, where $y=1$ corresponds to the star $S_n$.

	Let $4\le n\le12$, and let $r_n$ be the rational number in Table~\ref{tab:small}. This table also lists, rounded, the value $\rho_1$ of $\rho$ for the extremal graph and the largest value $\rho_2$ among the other graphs compared in the proofs. Since $r_n>n-2\ge x-1$, Lemma~\ref{lem:book} applies at $r_n$ to every $B_0[x,y]$ of order $n$. Exact rational evaluation shows that $F^{(0)}_{n,xy}(r_n)<0$ if $B_0[x,y]$ is $S_n$ for $n\le7$ or $H_{n,0}$ for $n\ge8$, and that $F^{(0)}_{n,xy}(r_n)>0$ for every other pair $(x,y)$. Thus, the graph named in Theorem~\ref{thm:main}(i) has $\rho>r_n$, and every other tree that can be a maximizer has $\rho<r_n$.

	Let $n\ge13$. By Lemma~\ref{lem:balance} with $q=0$ and $s=n\ge9$, every double star $B_0[x,y]$ with $x\ge y\ge2$ other than $H_{n,0}$ has a smaller $\rho$ than $H_{n,0}$. For the star, we have $\lfloor n^2/4\rfloor^2\ge(n^2-1)^2/16>(n-1)^3$, because $(n+1)^2>16(n-1)$ for $n\ge13$. Together with \eqref{eq:pf} applied to the spine edge of $H_{n,0}$, this gives
	\[
	\rho(S_n)=(n-1)^{3/2}<\left\lfloor\frac{n^2}4\right\rfloor<\rho(H_{n,0}).
	\]
	Hence, $H_{n,0}$ is the unique maximizer for $n\ge13$.
\end{proof}

\begin{table}[ht]
	\centering
	\caption{Separating values $r_n$ for small orders.}\label{tab:small}
	\small
	\setlength{\tabcolsep}{4pt}
	\begin{tabular}{rrrr@{\quad}rrrr@{\quad}rrrr}
		\toprule
		\multicolumn{4}{c}{Trees} & \multicolumn{4}{c}{Unicyclic graphs} & \multicolumn{4}{c}{Bicyclic graphs}\\
		$n$ & $r_n$ & $\rho_1$ & $\rho_2$ & $n$ & $r_n$ & $\rho_1$ & $\rho_2$ & $n$ & $r_n$ & $\rho_1$ & $\rho_2$\\
		\midrule
		4 & $5$ & 5.1962 & 4.8284 & 5 & $14.562$ & 14.5623 & 14.4557 & 5 & $21.5$ & 21.5506 & 20.5554\\
		5 & $7.8$ & 8.0000 & 7.5320 & 6 & $18.486$ & 18.4868 & 18.4781 & 6 & $26.3$ & 26.3303 & 26.0726\\
		6 & $11$ & 11.1803 & 10.6921 & 7 & $22.966$ & 22.9666 & 22.8326 & 7 & $31.4$ & 31.4191 & 30.8522\\
		7 & $14.5$ & 14.6969 & 14.3458 & 8 & $27.75$ & 27.7968 & 27.6830 & & & & \\
		8 & $18.55$ & 18.5830 & 18.5203 & 9 & $33$ & 33.2666 & 32.9706 & & & & \\
		9 & $23$ & 23.2184 & 22.7094 & 10 & $38.75$ & 39.0361 & 38.6530 & & & & \\
		10 & $28.3$ & 28.5078 & 28.1621 & 11 & $45$ & 45.4955 & 44.8141 & & & & \\
		11 & $34$ & 34.1254 & 33.3665 & 12 & $51.75$ & 52.2208 & 51.4575 & & & & \\
		12 & $40.2$ & 40.4499 & 40.0079 & & & & & & & & \\
		\bottomrule
	\end{tabular}
\end{table}

\section{Unicyclic graphs}\label{sec:uni}
For integers $x,y,z\ge2$, let $T[x,y,z]=T(x-2,y-2,z-2)$, so that the cycle vertices of $T[x,y,z]$ have degrees $x$, $y$ and $z$. Define the symmetric polynomial
\begin{equation}\label{eq:P}
	\begin{aligned}
		P_\lambda(x,y,z)
		&=\big(\lambda^2-(x-2)x^2\big)\big(\lambda^2-(y-2)y^2\big)\big(\lambda^2-(z-2)z^2\big)\\
		&\quad-\lambda^2\sum_{\mathrm{cyc}}x^2y^2\big(\lambda^2-(z-2)z^2\big)
		-2\lambda^3x^2y^2z^2.
	\end{aligned}
\end{equation}
Expanding the $3\times3$ determinant shows that $P_\lambda(x,y,z)=\det M_{T[x,y,z]}(\lambda)$. Since $K_3$ is complete multipartite with three parts of size one, Lemma~\ref{lem:det} gives the following criterion.

\begin{lemma}\label{lem:secular}
	Let $x,y,z\ge2$ and $\lambda>\max\{x,y,z\}-2$. Then $P_\lambda(x,y,z)$ has the sign of $\lambda-\rho(T[x,y,z])$.
\end{lemma}

Since $T[x,y,2]=B_1[x,y]$, Lemma~\ref{lem:book} gives $P_\lambda(x,y,2)=\lambda^2F^{(1)}_{x+y,xy}(\lambda)$. We now treat the graphs with pendent vertices at all three cycle vertices. The argument fixes the smallest cycle degree and uses concavity in the product of the other two.

\begin{lemma}\label{lem:concavity}
	Fix $z\ge2$ and $s=x+y$, where $x,y\ge2$, and set $S=s+z$. For $\lambda\ge S$, $P_\lambda(x,y,z)$ is a strictly concave function of $p=xy$ on the interval $2s-4\le p\le s^2/4$.
\end{lemma}

\begin{proof}
	Using $x^2+y^2=s^2-2p$ and $x^3+y^3=s^3-3sp$ in \eqref{eq:P}, we obtain a cubic in $p$. Differentiation yields
	\[
	-\frac12\frac{\partial^2P_\lambda}{\partial p^2}
	=\lambda^2\Big[\lambda^2+(2\lambda-S+6)z^2-3p+2s-4\Big]
	+(z-2)z^2(3p-2s+4).
	\]
	Since $p\le s^2/4$ and $\lambda\ge S>s$, we have $\lambda^2-3p>0$. Moreover, $2\lambda-S+6>0$, $2s-4>0$ because $s\ge4$, and $p\ge2s-4$ implies $3p-2s+4>0$. Hence, the right side is positive.
\end{proof}

For $x=y=t$ and $z=u$, the polynomial $P_\lambda$ factors as
\begin{equation}\label{eq:factor}
	P_\lambda(t,t,u)=\big[\lambda^2+t^2\lambda-(t-2)t^2\big]Q_\lambda(t,u),
\end{equation}
where
\begin{equation}\label{eq:Q}
	Q_\lambda(t,u)=
	\big[\lambda^2-t^2\lambda-(t-2)t^2\big]
	\big[\lambda^2-(u-2)u^2\big]-2t^2u^2\lambda^2.
\end{equation}
The first factor in \eqref{eq:factor} corresponds to the eigenvector of $M_{T[t,t,u]}(\lambda)$ that takes opposite values at the two cycle vertices of degree $t$, and it is positive for $\lambda\ge t$.

\begin{lemma}\label{lem:certU}
	Let $t,u\ge3$ be real numbers with $2t+u\ge16$. Then $Q_\lambda(t,u)>0$ for $\lambda=\lambda_{12}(2t+u-2)$.
\end{lemma}

\begin{proof}
	Put $R_U(t,u)=256\,Q_{\lambda_{12}(2t+u-2)}(t,u)$. This is a polynomial with integer coefficients, obtained by expanding \eqref{eq:Q} with integer additions and multiplications only. Each of the translated polynomials
	\[
	R_U(X+3,Y+6),\quad R_U(X+4,Y+4),\quad R_U(X+6,Y+3)
	\]
	has $44$ nonzero coefficients, all of them positive. Their constant coefficients are $8771328$, $405760$ and $14265153$, and the full coefficient arrays are listed in Tables~\ref{tab:RU1}--\ref{tab:RU3}. Hence, $R_U(t,u)>0$ if $t\ge3$ and $u\ge6$, if $t\ge4$ and $u\ge4$, or if $t\ge6$ and $u\ge3$. If $t<4$, then $u\ge16-2t>8$, and if $4\le t<6$, then $u\ge16-2t>4$. Therefore, these three regions cover all pairs in the lemma.
\end{proof}

\begin{lemma}\label{lem:interiorU}
	Let $n\ge13$, and let $a,b,c\ge1$ with $a+b+c=n-3$. Then $\rho(T(a,b,c))<\lambda_{12}(n+1)$.
\end{lemma}

\begin{proof}
	Let $x\ge y\ge z\ge3$ be the cycle degrees of $T(a,b,c)$, so that $x+y+z=S=n+3\ge16$, and put $\lambda=\lambda_{12}(n+1)=(S^2-2S+12)/4$. Then $\lambda>S$, since $S^2-6S+12>0$. Fix $z$ and put $s=x+y$. Since $z\le y\le s/2$, the product $p=xy$ satisfies $z(s-z)\le p\le s^2/4$. This interval lies in the interval of Lemma~\ref{lem:concavity}, because $z(s-z)-(2s-4)=(z-2)(s-z-2)\ge0$. By concavity, $P_\lambda(x,y,z)$ is at least the smaller of its values at the two endpoints, which correspond to the triples
	\[
	(S-2z,z,z),\quad ((S-z)/2,(S-z)/2,z).
	\]
	Up to permutation, both have the form $(t,t,u)$ with real $t,u\ge3$ and $2t+u=S$. By \eqref{eq:factor} and Lemma~\ref{lem:certU}, $P_\lambda$ is positive at both endpoints, because $\lambda>t$ and $2t+u-2=n+1$. Hence, $P_\lambda(x,y,z)>0$, and Lemma~\ref{lem:secular} gives $\rho(T(a,b,c))<\lambda$.
\end{proof}

\begin{proof}[Proof of Theorem~\ref{thm:main}(ii)]
	For $n=3$, the only unicyclic graph is $C_3=H_{3,1}$. Let $n\ge4$ and let $G$ be a graph with the largest $\rho$ in $\cU_n$. By Corollary~\ref{cor:families}, $G\cong T(a,b,c)$ with $a\ge b\ge c\ge0$ and $a+b+c=n-3$. For $n=4$, this gives $G\cong T(1,0,0)=H_{4,1}$.

	Let $n\ge13$. If $c=0$, then $G=B_1[a+2,b+2]$ with $b+2\ge2$, and Lemma~\ref{lem:balance} with $q=1$ and $s=n+1\ge12$ gives $G\cong H_{n,1}$. If $c\ge1$, then Lemmas~\ref{lem:interiorU} and~\ref{lem:lower} give $\rho(G)<\lambda_{12}(n+1)<\rho(H_{n,1})$, which is impossible. Hence, $G\cong H_{n,1}$.

	Let $8\le n\le12$ and $\lambda=\lambda_{12}(n+1)$. By Lemma~\ref{lem:lower}, $\rho(H_{n,1})>\lambda$. Exact rational evaluation shows that $P_\lambda(x,y,z)>0$ for every triple $x\ge y\ge z\ge2$ with $x+y+z=n+3$ other than the triple of $H_{n,1}$, except the four triples with $n=11$ and $z=2$. These are $33$ triples, and Lemma~\ref{lem:secular} applies since $\lambda>n+3$. The four excepted triples correspond to graphs $B_1[x,y]$ with $x+y=12$, which are covered by Lemma~\ref{lem:balance} with $s=12$. Accordingly, for $n=11$ the value $\rho_2$ in Table~\ref{tab:small} is taken over the triples with $z\ge3$. Hence, $G\cong H_{n,1}$.

	Let $5\le n\le7$, and let $r_n$ be the number in Table~\ref{tab:small}. Then $r_n>n+3$, and exact rational evaluation shows that $P_{r_n}(x,y,z)<0$ for the triple of $G_n^*$ and $P_{r_n}(x,y,z)>0$ for the other $6$ triples. By Lemma~\ref{lem:secular}, $\rho(G_n^*)>r_n>\rho(T[x,y,z])$ for every other triple, and hence $G\cong G_n^*$. In all cases, the maximizer is unique, which proves Theorem~\ref{thm:main}(ii), including the exceptional graph $T(1,1,1)$ for $n=6$.
\end{proof}

\section{Bicyclic graphs}\label{sec:bi}
For integers $x,y\ge3$ and $z\ge2$, let $D[x,y,z]=D(x-3,y-3,z-2)$, and for integers $x\ge3$ and $z\ge2$, let $K[x,z]=K(x-3,z-2)$. In the notation of Figure~\ref{fig:graphs}, the vertices $u,v,w_1,w_2$ of $D[x,y,z]$ have degrees $x,y,z,2$, and the vertices $u,v,w_1$ of $K[x,z]$ have degrees $x,3,z$, while the two remaining vertices of $K_{2,3}$ have degree $2$. The orders of $D[x,y,z]$ and $K[x,z]$ are $x+y+z-4$ and $x+z$, and $D[x,y,2]=B_2[x,y]$. Define
\begin{align}
	P^D_\lambda(x,y,z)&=\det\begin{pmatrix}
		\lambda^2-(x+1)x^2&-(\lambda+4)xy&-\lambda xz\\
		-(\lambda+4)xy&\lambda^2-(y+1)y^2&-\lambda yz\\
		-\lambda xz&-\lambda yz&\lambda^2-(z-2)z^2
	\end{pmatrix},\label{eq:PD}\\
	P^K_\lambda(x,z)&=\det\begin{pmatrix}
		\lambda^2-(x+5)x^2&-24x&-\lambda xz\\
		-24x&\lambda^2-72&-3\lambda z\\
		-\lambda xz&-3\lambda z&\lambda^2-(z-2)z^2
	\end{pmatrix}.\label{eq:PK}
\end{align}

\begin{lemma}\label{lem:detbi}
	The following statements hold.
	\begin{itemize}
		\item[(a)] $\det M_{D[x,y,z]}(\lambda)=\lambda^2P^D_\lambda(x,y,z)$. If $\lambda>\max\{x,y\}-3$ and $\lambda^2>(z-2)z^2$, then $P^D_\lambda(x,y,z)$ has the sign of $\lambda-\rho(D[x,y,z])$.
		\item[(b)] $\det M_{K[x,z]}(\lambda)=\lambda^4P^K_\lambda(x,z)$. If $\lambda^2>(x-3)x^2$ and $\lambda^2>(z-2)z^2$, then $P^K_\lambda(x,z)$ has the sign of $\lambda-\rho(K[x,z])$.
	\end{itemize}
\end{lemma}

\begin{proof}
	The graph $K_4-e$ is complete multipartite with parts $\{u\}$, $\{v\}$ and $\{w_1,w_2\}$, and $K_{2,3}$ is complete multipartite with parts $\{u,v\}$ and $\{w_1,w_2,w_3\}$. The stated conditions are those of Lemma~\ref{lem:det}. In $M_{D[x,y,z]}(\lambda)$, the row of $w_2$ is $(-2\lambda x,-2\lambda y,0,\lambda^2)$. Taking the Schur complement of its diagonal entry $\lambda^2$ produces the matrix in \eqref{eq:PD}, which gives (a). In $M_{K[x,z]}(\lambda)$, each of $w_2$ and $w_3$ has diagonal entry $\lambda^2$ and entries $-2\lambda x$ and $-6\lambda$ in the columns of $u$ and $v$. Taking the Schur complement of the corresponding block $\lambda^2I_2$ produces the matrix in \eqref{eq:PK}, which gives (b).
\end{proof}

Since $P^D_\lambda(x,y,z)$ is symmetric in $x$ and $y$, it is a polynomial in $s=x+y$, $p=xy$ and $z$. Expanding \eqref{eq:PD}, we obtain
\begin{equation}\label{eq:PDsp}
	\begin{aligned}
		P^D_\lambda={}&\lambda^6-\big[s^3+s^2-(3s+2)p+p^2+z^2(s^2-2p+z-2)\big]\lambda^4-2p^2(z^2+4)\lambda^3\\
		&+\big[p^3+(s-15)p^2+(s+z-8)p^2z^2+(z-2)z^2\big(s^3+s^2-(3s+2)p\big)\big]\lambda^2\\
		&+8(z-2)p^2z^2\lambda-(z-2)(p+s-15)p^2z^2.
	\end{aligned}
\end{equation}

\begin{lemma}\label{lem:concavityD}
	Fix integers $z\ge3$ and $s\ge6$, and put $n=s+z-4$. If $\lambda>n$, $\lambda^2>(z-2)z^2$ and $\lambda^2>3s^2/4+s$, then $P^D_\lambda$ is a strictly concave function of $p$ on the interval $3(s-3)\le p\le s^2/4$.
\end{lemma}

\begin{proof}
	Differentiating \eqref{eq:PDsp} twice with respect to $p$ gives
	\begin{align*}
		\frac{\partial^2P^D_\lambda}{\partial p^2}={}&-2z^2\lambda^2(2\lambda-s-z+8)-16\lambda\big(\lambda^2-(z-2)z^2\big)\\
		&-2\lambda^2\big(\lambda^2-3p-s+15\big)-2(z-2)z^2(3p+s-15).
	\end{align*}
	We have $2\lambda-s-z+8=2\lambda-n+4>0$ and $\lambda^2-(z-2)z^2>0$. Since $p\le s^2/4$, we have $\lambda^2-3p-s+15>\lambda^2-3s^2/4-s>0$. Finally, $z>2$ and $3p+s-15\ge10s-42>0$. Hence, all four terms are negative.
\end{proof}

For $x=y=t$, the polynomial $P^D_\lambda$ factors as
\begin{equation}\label{eq:factorD}
	P^D_\lambda(t,t,z)=\big[\lambda^2+t^2\lambda-(t-3)t^2\big]Q^D_\lambda(t,z),
\end{equation}
where
\begin{equation}\label{eq:QD}
	Q^D_\lambda(t,z)=\big[\lambda^2-t^2\lambda-(t+5)t^2\big]\big[\lambda^2-(z-2)z^2\big]-2t^2z^2\lambda^2.
\end{equation}
As in \eqref{eq:factor}, the first factor in \eqref{eq:factorD} is positive for $\lambda\ge t$. We use the following three polynomials with integer coefficients, obtained by substituting the value $\lambda_{27}(n+2)$ for the corresponding order $n$:
\begin{align*}
	R_1(t,z)&=256\,Q^D_{\lambda_{27}(2t+z-2)}(t,z),\\
	R_2(x,z)&=4096\,P^D_{\lambda_{27}(x+z+1)}(x,3,z),\\
	R_3(x,z)&=4096\,P^K_{\lambda_{27}(x+z+2)}(x,z).
\end{align*}

\begin{lemma}\label{lem:certB}
	The following statements hold.
	\begin{itemize}
		\item[(a)] $R_1(t,z)>0$ for all real $t,z\ge3$ with $2t+z\ge12$.
		\item[(b)] $R_2(x,z)>0$ for all real $x,z\ge3$.
		\item[(c)] $R_3(x,z)>0$ for all real $x\ge3$ and all $z\in\{2\}\cup[3,\infty)$.
	\end{itemize}
\end{lemma}

\begin{proof}
	For (a), each of $R_1(X+4,Y+3)$ and $R_1(X+3,Y+4)$ has $44$ nonzero coefficients, all of them positive, with constant coefficients $13495248$ and $17731973$. The arrays are listed in Tables~\ref{tab:R11} and~\ref{tab:R12}. Hence, $R_1(t,z)>0$ if $t\ge4$ and $z\ge3$, or if $t\ge3$ and $z\ge4$. If $t<4$ and $z<4$, then $2t+z<12$, so these two regions cover all pairs in (a).

	For (b), the leading form of $R_2$ vanishes on the line $x=z$, so translations alone cannot produce positive coefficients, and we use two substitutions adapted to this line. Each of $R_2(Y+3+X,Y+3)$ and $R_2(X+3,X+3+Y)$ has $89$ nonzero coefficients, all of them positive, and both constant coefficients equal $96454048320$. The first polynomial gives $R_2(x,z)>0$ for $x\ge z\ge3$, with $X=x-z$ and $Y=z-3$, and the second gives $R_2(x,z)>0$ for $z\ge x\ge3$.

	For (c), the same substitutions give polynomials $R_3(Y+3+X,Y+3)$ and $R_3(X+3,X+3+Y)$ with $89$ positive coefficients each and constant coefficient $843887666233$, which proves (c) for $z\ge3$. Finally, $R_3(X+3,2)$ is a polynomial in $X$ of degree $12$ with $13$ positive coefficients, which proves (c) for $z=2$.

	Table~\ref{tab:certsum} summarizes these data.
\end{proof}

\begin{lemma}\label{lem:interiorB}
	Let $n\ge5$ and $\lambda=\lambda_{27}(n+2)$.
	\begin{itemize}
		\item[(a)] If $n\ge8$, then $\rho(D[x,y,z])<\lambda$ for every graph $D[x,y,z]$ of order $n$ with $z\ge3$.
		\item[(b)] $\rho(K[x,z])<\lambda$ for every graph $K[x,z]$ of order $n$.
	\end{itemize}
\end{lemma}

\begin{proof}
	We first collect three inequalities. Since $\lambda>(n+2)^2/4$, we have $\lambda>n$, and $\lambda^2>(n+2)^4/16>n^3$, because the minimum of $(n+2)^4/n^3$ over $n>0$ is $4096/216>16$, attained at $n=6$. Moreover, $(n+2)^4/16\ge(n+2)^2>3(n+1)^2/4+n+1$ for $n\ge2$.

	(a) By symmetry, we may assume that $x\ge y$. Fix $z\in[3,n-2]$ and put $s=x+y=n+4-z\ge6$. Then $(z-2)z^2<n^3<\lambda^2$, $s\le n+1$ and $\lambda^2>3s^2/4+s$, so Lemma~\ref{lem:concavityD} applies. The product $p=xy$ lies in $[3(s-3),s^2/4]$, and by concavity, $P^D_\lambda(x,y,z)$ is at least the smaller of its values at $p=3(s-3)$ and $p=s^2/4$. The first endpoint corresponds to $(s-3,3,z)$, where $P^D_\lambda(s-3,3,z)=R_2(s-3,z)/4096>0$ by Lemma~\ref{lem:certB}(b), since $(s-3)+z+1=n+2$. The second endpoint corresponds to $(s/2,s/2,z)$. By \eqref{eq:factorD}, $P^D_\lambda(s/2,s/2,z)$ is the product of a positive factor and $Q^D_\lambda(s/2,z)=R_1(s/2,z)/256$, which is positive by Lemma~\ref{lem:certB}(a), since $s+z=n+4\ge12$. Hence, $P^D_\lambda(x,y,z)>0$. Since $\lambda>n>\max\{x,y\}-3$, Lemma~\ref{lem:detbi}(a) gives $\rho(D[x,y,z])<\lambda$.

	(b) We have $x\le n-2$ and $z\le n-3$, so $(x-3)x^2<\lambda^2$ and $(z-2)z^2<\lambda^2$. Since $x+z=n$, Lemma~\ref{lem:certB}(c) gives $P^K_\lambda(x,z)=R_3(x,z)/4096>0$, and Lemma~\ref{lem:detbi}(b) gives $\rho(K[x,z])<\lambda$.
\end{proof}

\begin{proof}[Proof of Theorem~\ref{thm:main}(iii)]
	For $n=4$, the only bicyclic graph is $K_4-e=H_{4,2}$. Let $n\ge5$, and let $G$ be a graph with the largest $\rho$ in $\cB_n$. By Corollary~\ref{cor:families}, $G$ is isomorphic to some $D[x,y,z]$ with $x\ge y\ge3$ and $z\ge2$, or to some $K[x,z]$ with $x\ge3$ and $z\ge2$.

	Let $n\ge8$. If $G\cong D[x,y,2]=B_2[x,y]$, then Lemma~\ref{lem:balance} with $q=2$ and $s=n+2\ge10$ gives $G\cong H_{n,2}$. If $G\cong D[x,y,z]$ with $z\ge3$ or $G\cong K[x,z]$, then Lemmas~\ref{lem:interiorB} and~\ref{lem:lower} give $\rho(G)<\lambda_{27}(n+2)<\rho(H_{n,2})$, which is impossible. Hence, $G\cong H_{n,2}$.

	Let $5\le n\le7$, and let $r_n$ be the number in Table~\ref{tab:small}. For every graph $D[x,y,z]$ and $K[x,z]$ of order $n$, the conditions of Lemma~\ref{lem:detbi} hold at $\lambda=r_n$. Exact rational evaluation shows that $P^D_{r_n}<0$ for $H_{n,2}$, and that $P^D_{r_n}>0$ and $P^K_{r_n}>0$ for the other $15$ graphs. Hence, $G\cong H_{n,2}$. In all cases, the maximizer is unique.
\end{proof}

\section{Concluding remarks}\label{sec:concl}
Theorem~\ref{thm:main} shows that, for $q\in\{0,1,2\}$ and all orders $n$ outside the small exceptions, the graph $H_{n,q}$ uniquely maximizes the second Zagreb spectral radius among connected graphs of order $n$ with cyclomatic number $q$. This agrees with a simple degree count. If $uv$ is an edge of a connected graph $G$ of order $n$, then $d_u+d_v\le n+|N(u)\cap N(v)|$. The subgraph formed by $u$, $v$ and their common neighbours has cyclomatic number $|N(u)\cap N(v)|$, so $d_u+d_v\le n+q$ when $G$ has cyclomatic number $q$. Equality forces every vertex other than $u$ and $v$ to be adjacent to $u$ or $v$, and forces $G$ to have no further edges, that is, $G$ is obtained from $K_2\vee\Kb$ by attaching pendent vertices to $u$ and $v$. Hence, $H_{n,q}$ maximizes the largest entry $d_ud_v$ of $Z(G)$, and by Theorem~\ref{thm:values}, $\rho(H_{n,q})$ exceeds this entry by about $(n+q)/2+7q-2$.

For $q=3$, the situation is different for small orders. Let $K_4^{(n)}$ be the graph obtained from $K_4$ by attaching $n-4$ pendent vertices to two of its vertices as evenly as possible. An exhaustive computation over all connected graphs of order $n\le12$ with cyclomatic number $3$ shows that $K_4^{(n)}$ is the unique maximizer of $\rho$ for $5\le n\le12$. A direct comparison of the two graphs shows that $\rho(H_{n,3})<\rho(K_4^{(n)})$ for $5\le n\le14$ and $\rho(H_{n,3})>\rho(K_4^{(n)})$ for $15\le n\le200$. These observations and Theorem~\ref{thm:main} suggest the following conjecture.

\begin{conjecture}\label{conj:c}
	For every integer $q\ge3$, there is an integer $n_q$ such that $H_{n,q}$ is the unique graph with the largest second Zagreb spectral radius among connected graphs of order $n\ge n_q$ with cyclomatic number $q$.
\end{conjecture}

The proof of Lemma~\ref{lem:det} uses only the product form $d_id_j$ of the weights. For weights $f(x,y)=g(x)g(y)$ with $g>0$, the identity \eqref{eq:mdl} and the sign criterion in Lemma~\ref{lem:det} remain valid after replacing $d_i$ by $g(d_i)$ in $D$ and $\bd$, and $p_id_i^2$ by $p_ig(1)^2g(d_i)^2$. The sufficient condition for $N_G(\lambda)$ to be positive definite then becomes $\lambda>p_ig(1)^2$ for parts $\{i\}$ of size one and $\lambda^2>p_ig(1)^2g(d_i)^2$ for the other vertices. For example, this applies to $f(x,y)=(xy)^\alpha$. For $\alpha\ge1$, this function is also increasing and convex in each variable, so Lemmas~\ref{lem:trees} and~\ref{lem:shen} apply as well. It would be interesting to determine the maximizers for these weights, and more generally for the weights with property $P^{**}$ considered in \cite{Shen2025}.

\appendix
\section{Certificate data}\label{app:coeff}
Table~\ref{tab:certsum} summarizes the certificates of Lemmas~\ref{lem:certU} and~\ref{lem:certB}, and Tables~\ref{tab:RU1}--\ref{tab:R12} give the complete expansions $\sum_{i,j}c_{ij}X^iY^j$ of the polynomials used in Lemmas~\ref{lem:certU} and~\ref{lem:certB}(a). All coefficients were obtained by exact integer arithmetic.

\begin{table}[H]
	\centering
	\caption{Summary of the certificates.}\label{tab:certsum}
	\begin{tabular}{llrr}
		\toprule
		Polynomial & Region certified & Nonzero coefficients & Constant coefficient\\
		\midrule
		$R_U(X+3,Y+6)$ & $t\ge3$, $u\ge6$ & 44 & 8771328\\
		$R_U(X+4,Y+4)$ & $t\ge4$, $u\ge4$ & 44 & 405760\\
		$R_U(X+6,Y+3)$ & $t\ge6$, $u\ge3$ & 44 & 14265153\\
		$R_1(X+4,Y+3)$ & $t\ge4$, $z\ge3$ & 44 & 13495248\\
		$R_1(X+3,Y+4)$ & $t\ge3$, $z\ge4$ & 44 & 17731973\\
		$R_2(Y+3+X,Y+3)$ & $x\ge z\ge3$ & 89 & 96454048320\\
		$R_2(X+3,X+3+Y)$ & $z\ge x\ge3$ & 89 & 96454048320\\
		$R_3(Y+3+X,Y+3)$ & $x\ge z\ge3$ & 89 & 843887666233\\
		$R_3(X+3,X+3+Y)$ & $z\ge x\ge3$ & 89 & 843887666233\\
		$R_3(X+3,2)$ & $x\ge3$, $z=2$ & 13 & 304692840000\\
		\bottomrule
	\end{tabular}
\end{table}

\begin{table}[H]
\centering
\setlength{\tabcolsep}{3pt}
\caption{Coefficients of $R_U(X+3,Y+6)$. A dash denotes zero.}\label{tab:RU1}
\small
\begin{tabular}{r|rrrrrrrrr}
\toprule
$i\backslash j$ & $0$ & $1$ & $2$ & $3$ & $4$ & $5$ & $6$ & $7$ & $8$ \\
\midrule
0 & 8771328 & 11332224 & 6474816 & 1948320 & 336960 & 34936 & 2148 & 72 & 1 \\
1 & 4810752 & 8941824 & 5510016 & 1585408 & 241152 & 20224 & 888 & 16 & - \\
2 & 463104 & 2094336 & 1428800 & 398432 & 51776 & 3192 & 76 & - & - \\
3 & 866304 & 577024 & 188160 & 41920 & 4160 & 144 & - & - & - \\
4 & 676608 & 279040 & 31936 & 2496 & 112 & - & - & - & - \\
5 & 186880 & 67840 & 5248 & 128 & - & - & - & - & - \\
6 & 22528 & 7040 & 320 & - & - & - & - & - & - \\
7 & 1024 & 256 & - & - & - & - & - & - & - \\
\bottomrule
\end{tabular}
\end{table}

\begin{table}[H]
\centering
\setlength{\tabcolsep}{3pt}
\caption{Coefficients of $R_U(X+4,Y+4)$. A dash denotes zero.}\label{tab:RU2}
\small
\begin{tabular}{r|rrrrrrrrr}
\toprule
$i\backslash j$ & $0$ & $1$ & $2$ & $3$ & $4$ & $5$ & $6$ & $7$ & $8$ \\
\midrule
0 & 405760 & 290304 & 1008768 & 792832 & 212400 & 28096 & 1992 & 72 & 1 \\
1 & 521216 & 759296 & 634880 & 581760 & 145152 & 15904 & 816 & 16 & - \\
2 & 2344704 & 2003328 & 286080 & 153824 & 33248 & 2712 & 76 & - & - \\
3 & 2028032 & 1524480 & 154752 & 22080 & 3168 & 144 & - & - & - \\
4 & 747776 & 514944 & 46848 & 2240 & 112 & - & - & - & - \\
5 & 140288 & 88320 & 6400 & 128 & - & - & - & - & - \\
6 & 13312 & 7552 & 320 & - & - & - & - & - & - \\
7 & 512 & 256 & - & - & - & - & - & - & - \\
\bottomrule
\end{tabular}
\end{table}

\begin{table}[H]
\centering
\setlength{\tabcolsep}{3pt}
\caption{Coefficients of $R_U(X+6,Y+3)$. A dash denotes zero.}\label{tab:RU3}
\small
\begin{tabular}{r|rrrrrrrrr}
\toprule
$i\backslash j$ & $0$ & $1$ & $2$ & $3$ & $4$ & $5$ & $6$ & $7$ & $8$ \\
\midrule
0 & 14265153 & 28378080 & 585756 & 780960 & 358662 & 50464 & 3228 & 96 & 1 \\
1 & 24455952 & 44245872 & 2727120 & 527920 & 193584 & 22096 & 1008 & 16 & - \\
2 & 17291412 & 27552048 & 2171876 & 164768 & 38204 & 3120 & 76 & - & - \\
3 & 6305424 & 8974864 & 724320 & 30304 & 3344 & 144 & - & - & - \\
4 & 1296432 & 1672960 & 120160 & 3072 & 112 & - & - & - & - \\
5 & 152704 & 180352 & 9856 & 128 & - & - & - & - & - \\
6 & 9664 & 10496 & 320 & - & - & - & - & - & - \\
7 & 256 & 256 & - & - & - & - & - & - & - \\
\bottomrule
\end{tabular}
\end{table}

\begin{table}[H]
\centering
\setlength{\tabcolsep}{3pt}
\caption{Coefficients of $R_1(X+4,Y+3)$. A dash denotes zero.}\label{tab:R11}
\small
\begin{tabular}{r|rrrrrrrrr}
\toprule
$i\backslash j$ & $0$ & $1$ & $2$ & $3$ & $4$ & $5$ & $6$ & $7$ & $8$ \\
\midrule
0 & 13495248 & 13595328 & 3656928 & 846528 & 159688 & 20720 & 1576 & 64 & 1 \\
1 & 5216256 & 10201728 & 2081536 & 488512 & 101376 & 12064 & 704 & 16 & - \\
2 & 173664 & 4616384 & 588048 & 105216 & 23176 & 2256 & 76 & - & - \\
3 & 36160 & 1756224 & 175296 & 14272 & 2448 & 144 & - & - & - \\
4 & 124864 & 475648 & 44352 & 1792 & 112 & - & - & - & - \\
5 & 40704 & 78080 & 6016 & 128 & - & - & - & - & - \\
6 & 5248 & 6912 & 320 & - & - & - & - & - & - \\
7 & 256 & 256 & - & - & - & - & - & - & - \\
\bottomrule
\end{tabular}
\end{table}

\begin{table}[H]
\centering
\setlength{\tabcolsep}{3pt}
\caption{Coefficients of $R_1(X+3,Y+4)$. A dash denotes zero.}\label{tab:R12}
\small
\begin{tabular}{r|rrrrrrrrr}
\toprule
$i\backslash j$ & $0$ & $1$ & $2$ & $3$ & $4$ & $5$ & $6$ & $7$ & $8$ \\
\midrule
0 & 17731973 & 12439936 & 3976520 & 895608 & 148962 & 17520 & 1312 & 56 & 1 \\
1 & 11004456 & 8484048 & 2685160 & 653488 & 110680 & 11632 & 664 & 16 & - \\
2 & 2362964 & 2356136 & 606932 & 157648 & 26764 & 2280 & 76 & - & - \\
3 & 415072 & 642320 & 90240 & 17824 & 2720 & 144 & - & - & - \\
4 & 191408 & 222016 & 23200 & 1600 & 112 & - & - & - & - \\
5 & 60800 & 50560 & 4480 & 128 & - & - & - & - & - \\
6 & 8896 & 5760 & 320 & - & - & - & - & - & - \\
7 & 512 & 256 & - & - & - & - & - & - & - \\
\bottomrule
\end{tabular}
\end{table}

\section*{Acknowledgment}
This research was supported by the Higher Education Scientific Research Project of the Anhui Provincial Department of Education under Grants 2025AHGXZK31124, 2025AHGXZK31313 and SK2017A0586.

\section*{Declarations}

\noindent\textbf{Data availability} All data generated or analysed during this study are included in this article.  \\

\noindent\textbf{Use of AI tools} The authors used the AI assistant Claude (Anthropic) to assist with parts of this work, including computations and the preparation of the manuscript. The authors have checked all content and take full responsibility for this article.  \\

\noindent\textbf{Conflict of Interest} The authors declare no conflict of interest.

\end{sloppypar}
\end{document}